\documentclass[12pt]{article}
\usepackage{amsfonts,amssymb}
\newtheorem{corollary}{Corollary}
\newtheorem{proposition}{Proposition}
\newtheorem{theorem}{Theorem}
\newtheorem{definition}{Definition}
\newtheorem{example}{Example}

\begin{document}
\begin{center}
{\large \bf DYNAMICS OF FINITE-MULTIVALUED TRANSFORMATIONS}
\vspace{0.5cm}

 KONSTANTIN IGUDESMAN
 \footnote{The author was supported in part by RF Education
Ministry and DAAD grant 331 4 00 088}
\vspace{0.2cm}

E-mail: Konstantin.Igudesman@ksu.ru
\end{center}
{\bf ABSTRACT}\\
We consider a transformation of a normalized measure space such
that the image of any point is a finite set. We call such
transformation $m$-transformation. In this case the orbit of any
point looks like a tree. In the study of $m$-transformations we
are interested in the properties of the trees.

An $m$-transformation generates a stochastic kernel and a new
measure. Using these objects, we introduce analogies of some main
concept of ergodic theory: ergodicity, Koopman and
Frobenius-Perron operators etc. We prove ergodic theorems and
consider examples. We also indicate possible applications to
fractal geometry and give a generalization of our construction.
Some results which have analogies in the classical ergodic theory
we are proved using standard methods (see \cite{Boy}, \cite{Las}).
Other results, for instance Theorem \ref{equivalent} and Example
\ref{ex5}, have no analogies.

\section{Main definitions and examples}

Throughout the paper $(X,\mathcal{B},\mu)$ denotes a normalized
measure space. Let $m$ be a positive integer.

\begin{definition}
We call a multivalued transformation $S:X\rightarrow X$ an
$\mathbf{m}$-{\bf transformation} if $1\leq|S(x)|\leq m$ for any
$x\in X$, where $|A|$ is just a number of elements in $A$.
\end{definition}

Let
$$
S^{-1}_{k;l}(B)\equiv \{x\in X : |S(x)|=k, |S(x)\cap B|=l\},
$$
where $B\subset X$ and $k,l\in \mathbb{N}$. Note that sets
$S^{-1}_{k;l}(B)$ are pairwise disjoint for the fixed $B$.

\begin{definition}\label{measurable}
The $m$-transformation $S:X\rightarrow X$ is {\bf measurable} if
$S^{-1}_{k;l}(B)\in \mathcal{B}$ for all $B\in \mathcal{B}$ and
$k,l\in \mathbb{N}$.
\end{definition}

Let $K:X\times \mathcal{B}\rightarrow \mathbb{R}^+$ be the
function
$$
K(x,B)\equiv \frac{1}{|S(x)|}\sum _{y\in S(x)}\chi _B(y)\ .
$$
For each $x\in X$, $K(x,\cdot):\mathcal{B}\rightarrow
\mathbb{R}^+$ is a normalized measure and for each $B\in
\mathcal{B}$, $K(\cdot ,B):X \rightarrow \mathbb{R}^+$ is
measurable by the Definition \ref{measurable}. Therefore $K$ is a
{\bf stochastic kernel} that describes the $m$-transformation $S$.
We will use $K$ as a tool for proving some results. Fore a more
complete study of stochastic kernels the reader is referred to
\cite{Kre}.

For any measurable $m$-transformation $S$ we define a new measure
$S\mu$ on $(X,\mathcal{B},\mu)$
$$
S\mu (B)\equiv \int \limits_X K(x,B)\ {\rm d}\mu =\sum
_{k=1}^m\sum _{l=1}^k \frac{l}{k}\ \mu (S^{-1}_{k;l}(B)).
$$

\begin{definition}
We say the measurable $m$-transformation $S:X\rightarrow X$ {\bf
preserves measure} $\mu $ or that $\mu $ is $S$-{\bf invariant} if
$S\mu =\mu$.
\end{definition}

\begin{definition}
Let the $m$-transformation $S:X\rightarrow X$ preserve measure
$\mu $. The quadruple $(X,\mathcal{B},\mu ,S)$ is called an {\bf
m-dynamical system}.
\end{definition}

The next proposition gives a number of examples of $m$-dynamical
systems.

\begin{proposition}
Let $\{ S_i \}_1^k$ be a finite collection of the $\mu$-preserving
$m_i$-transformations of $(X,\mathcal{B},\mu)$ and let
$S(x)=\bigcup _{i=1}^kS_i(x)$ be measurable. Let $K, K_i$ be the
stochastic kernels that generates $S, S_i$ correspondently. If for
any $B\in \mathcal{B}$
\begin{equation}\label{union}
 K(x,B)=\frac{1}{k}\sum _{i=1}^kK_i(x,B)
\end{equation}
 for almost all $x\in X$, then $S$ is $\mu$-preserving.
\end{proposition}

$\blacktriangleright $ For any measurable $B$ we have
$$
S\mu (B)=\int \limits_X K(x,B)\ {\rm d}\mu =\frac{1}{k}\sum
_{i=1}^k \int \limits_X P_i(x,B)\ {\rm d}\mu =\mu (B)\ . \
\blacktriangleleft
$$

In the following examples $\lambda$ denotes the Lebesgue measure
on $[0,1]$.

\begin{example}\label{ex1}
Let $S:[0,1]\rightarrow [0,1]$ be defined by $S(x)=\{ x,1-x \}$.
Then $S$ is $\lambda $-preserving.
\end{example}

\begin{example}\label{ex2}
Let $S:[0,1]\rightarrow [0,1]$ be defined by
$$
S(x)=\left\{
\begin{array}{cc}
  \{ 2x,1-2x \}\ , & x\in [0,\frac{1}{2}]\phantom{1}  \\
   & \\
  \{ 2x-1 \}\ , & x\in (\frac{1}{2},1]\ .
\end{array}
\right.
$$
 Then $S$ is $\lambda $-preserving.
\end{example}

The following example shows that not every $\lambda$-preserving
$m$-transformation is union of $\lambda$-preserving
transformations.

\begin{example}\label{ex3}
Let $S:[0,1]\rightarrow [0,1]$ be defined by
$$
S(x)=\left\{
\begin{array}{cc}
  \{ \frac{3}{2}x \}\ , & x\in [0,\frac{1}{3})\phantom{1}  \\
   & \\
  \{ \frac{3}{2}x, \frac{3}{2}x-\frac{1}{2} \}\ , & x\in
  [\frac{1}{3},\frac{2}{3}] \\ & \\
  \{ \frac{3}{2}x-\frac{1}{2} \}\ , & x\in (\frac{2}{3},1]\ .
\end{array}
\right.
$$
Then $S$ is $\lambda $-preserving, but $S$ can not be represented
as union of $\lambda $-preserving transformations.
\end{example}

$\blacktriangleright $\ Assume $S(x)=\cup _{i=1}^kS_i(x)$, where
$S_i$ are the $\lambda $-preserving transformations. Then there
are a measurable set $B\subset [\frac{1}{3},\frac{2}{3}]$ of
positive measure and transformation $S_i$ (for instance $S_1$),
such that $S_1(B)\subset [0,\frac{1}{2}]$. We have
$$
\lambda (S_1^{-1}(S_1(B)))=\lambda (B\cup
(B-\frac{1}{3}))=2\lambda (B)\ \mbox{and}\  \lambda
(S_1(B))=\frac{3}{2}\lambda (B)\ .
$$
Since $S_1$ is the $\lambda$-preserving transformation, $\lambda
(S_1(B))=\lambda (B)=0$. \ $\blacktriangleleft $

\begin{example}
Let $S:[0,1]\rightarrow [0,1]$ be defined by
$$
S(x)=\left\{
\begin{array}{cc}
  \{ 2x,1-2x,x \}\ , & x\in [0,\frac{1}{2}]\phantom{1}  \\
   & \\
  \{ 2x-1,x \}\ , & x\in (\frac{1}{2},1]\ .
\end{array}
\right.
$$
 Then $S$ isn't $\lambda $-preserving.
\end{example}

$\blacktriangleright $\ For instance,
$$
S\lambda ([0,\frac{1}{2}])=\frac{2}{3}\lambda
([0,\frac{1}{2}])+\frac{1}{2}\lambda
([\frac{1}{2},\frac{3}{4}])=\frac{11}{24}\neq \lambda
([0,\frac{1}{2}])\ .
$$
Nevertheless, we can represent $S$ as the union of the
$\lambda$-preserving transformations $S_1(x)=x$ and $S_2$ from
Example \ref{ex2}. Of course, (\ref{union}) does not hold true.\
$\blacktriangleleft $

Let $S^{-1}(B)=\{ x\in X : S(x)\cap B\neq \emptyset \}$ denote the
full preimage of $B$.

\begin{definition}
A measurable $m$-transformation $S:X\rightarrow X$ is said to be
{\bf nonsingular} if for any $B\in \mathcal{B}$ such that $\mu
(B)=0$, we have $\mu (S^{-1}(B))=0$, i.e., $S\mu \ll \mu$.
\end{definition}

\section{Recurrence and ergodic theorems}

Let $S:X\rightarrow X$ be an $m$-transformation. The $n$-th
iterate of $S$ is denoted by $S^n$. The {\bf tree} at $x_0\in X$
is the set $\{x\in X : x\in S^n(x_0)\ \mbox{for some}\ n\geq 0\}$.
Any sequence $x_0, x_1, x_2, \ldots $ with $x_{n+1}\in S(x_n)\
\mbox{for all}\  n\geq 0$ is called {\bf orbit} of $x_0$.

In the study of $m$-dynamical systems, we are interested in
properties of the trees. For example, in the recurrence of trees
of $S$, i.e., the property that if the tree in $x$ starts in a
specified set, some orbits of $x$ return to that set infinitely
many times.

\begin{proposition}
Let $S$ be a nonsingular $m$-transformation on
$(X,\mathcal{B},\mu)$ and let $\mu (A)\leq \mu (S^{-1}(A))$ for
any $A\in \mathcal{B}$. If $\mu (B)>0$, then for almost all $x\in
B$ there is an orbit of $x$ that returns infinitely often to $B$.
\end{proposition}

$\blacktriangleright $\ Let $B$ be a measurable set with $\mu
(B)>0$, and let us define the set $A$ of points that never return
to $B$, i.e., $A=\{ x\in B : S^n(x)\cap B =\emptyset \ \mbox{for
all}\ n\geq 1\}=B\backslash \cup _{n=1}^\infty S^{-n}(B)$.
Consider a collection of sets
$$
A_1=A\cup S^{-1}(A),\ A_i=A\cup S^{-1}(A_{i-1}),\ i\geq 2\ .
$$
It is clear that $A\cap S^{-1}(A_{i-1})=\emptyset $. Hence
$$
\mu (A_i)=\mu (A)+\mu(S^{-1}(A_{i-1}))\geq \mu
(A)+\mu(A_{i-1})\geq \ldots \geq(i+1)\mu (A)\ .
$$
Therefore, $\mu (A)=0$. Since $\mu $ is nonsingular, $\mu
(S^{-n}(A))=0$ for any $n\geq 0$. This gives $\mu (B\backslash
\bigcup _n S^{-n}(A))=\mu (B)$, and for any $x\in B\backslash
\bigcup _n S^{-n}(A)$ there exists a orbit of $x$ that return
infinitely often to $B$. \ $\blacktriangleleft $

If $S$ is measure preserving, then we have an analogue of
Poincare's Recurrence Theorem.

\begin{corollary}\label{Poincare}
Let $S$ be a measure-preserving $m$-transformation on
$(X,\mathcal{B},\mu)$. If $\mu (B)>0$, then for almost all $x\in
B$ there is an orbit of $x$ that returns infinitely often to $B$.
\end{corollary}

$\blacktriangleright $\ Note that $S\mu \ll \mu$ and for any
measurable $A$
$$
\mu (A)=S\mu (A)=\sum _{k=1}^m\sum _{l=1}^k \frac{l}{k}\ \mu
(S^{-1}_{k;l}(A))\leq \mu (S^{-1}(A))\ .\ \blacktriangleleft
$$

Example \ref{ex1} shows there are orbits that do not return to
$B$. If $B=[0,\frac{1}{2})$, then for any $x\in B$ the orbit
$\{x,1-x,1-x,\ldots\}$ doesn't return to $B$.

For any nonsingular $m$-transformation $S$ and function $f$ on $X$
we define a new function $Uf$ on $X$ by the equality
$$
Uf(x)\equiv \int \limits_X f\ {\rm d}K(x,\cdot )
=\frac{1}{|S(x)|}\sum _{y\in S(x)}f(y)\ .
$$

\begin{proposition}
If $S$ is a nonsingular $m$-transformation and $f$ is a
real-valued measurable function on $X$, then
$$
\int \limits_X f\ {\rm d}S\mu=\int \limits_X Uf\ {\rm d}\mu \ ,
$$
in the sense that if one of these integrals exists then so does
the other and the two are equal.
\end{proposition}

$\blacktriangleright $\ We first show that $Uf$ is measurable.
Given any $\alpha \in \mathbb{R}$ consider an increasing sequence
of rational numbers $\alpha _1< \ldots <\alpha _k$, where $k\leq
m$ and $\sum _{i=1}^k \alpha _i<k \alpha$. Then the set
$$
B_{\alpha _1, \ldots ,\alpha _k}=S^{-1}(f^{-1}(-\infty,\alpha
_1])\cap S^{-1}(f^{-1}(\alpha _1,\alpha _2])\cap \ldots \cap
S^{-1}(f^{-1}(\alpha _{k-1},\alpha _k])
$$
is measurable. Taking the union of $B_{\alpha _1, \ldots ,\alpha
_k}$ for all possible $k\leq m$ and $\alpha _1, \ldots ,\alpha _k$
we conclude that the set $\{ x: (Uf)(x)<\alpha \}$ is measurable.

When $f=\chi _B$ is the characteristic function of $B\in
\mathcal{B}$,
$$
\int \limits_X \chi _B\ {\rm d}S\mu=S\mu (B)
$$
and
$$
\int \limits_X U\chi _B\ {\rm d}\mu =\int \limits_X
\frac{1}{|S(x)|}\sum _{y\in S(x)}\chi _B(y)\ {\rm d}\mu =\int
\limits_X\sum _{k=1}^m\sum _{l=0}^k \frac{l}{k}\ \chi
_{(S^{-1}_{k;l}(B))}\ {\rm d}\mu =S\mu (B)\ .
$$
Since $U$ is a linear operator, the formula is also true for
simple functions. If $f$ is a nonnegative measurable function,
then $f$ is the $S\mu$-pointwise limit of an increasing sequence
of simple functions $f_i$, and the result follows from the fact
that $Uf$ is the $\mu$-pointwise limit of the increasing sequence
of functions $Uf_i$ and Monotone Converges Theorem. Finally, any
measurable function $f$ can be written as the difference
$f=f^+-f^-$ of two nonnegative measurable functions, so the
formula is true in general.  \ $\blacktriangleleft $

\begin{corollary}\
Let $S:X\rightarrow X$ be a measurable $m$-transformation on
$(X,\mathcal{B},\mu)$. Then $S$ is $\mu$-preserving if and only if
$$
\int \limits_X f\ {\rm d}\mu=\int \limits_X Uf\ {\rm d}\mu
$$
for any $f\in \mathcal{L}^1$.
\end{corollary}

$\blacktriangleright $\ This follows from the Proposition above
and from the equality
$$
\mu (B)=\int \limits_X U\chi _B\ {\rm d}\mu =\int \limits_X \left(
\int \limits_X \chi _B\ {\rm d}K(x,\cdot )\right)\ {\rm d}\mu
=\int \limits_X K(x,B)\ {\rm d}\mu =S\mu(B)\ .\ \blacktriangleleft
$$

\begin{proposition}
Let $S:X\rightarrow X$ be a $\mu$-preserving $m$-transformation on
$(X,\mathcal{B},\mu)$. Then the positive linear operator $U$ is a
contraction on $\mathcal{L}^p$ for any $1\leq p\leq \infty$.
\end{proposition}

$\blacktriangleright $\ It is easily seen that $U$ is a
contraction on $\mathcal{L}^{\infty}$. By Jensen inequality
$|Uf|^p\leq U|f|^p$ for any $p\geq 1$ and $f\in \mathcal{L}^p$
(see \cite{Kre}, Chapter 1, Lemma 7.4 for a more general
statement). Then
$$
\| Uf \| _p^p=\int \limits_X |Uf|^p\ {\rm d}\mu \leq \int
\limits_X U|f|^p\ {\rm d}\mu =\int \limits_X |f|^p\ {\rm d}\mu =\|
f \| _p^p \ .\ \blacktriangleleft
$$

For a function $f$ on $X$ and an $m$-transformation
$S:X\rightarrow X $, we define the averages
$$
A_n(f)=\frac{1}{n}\sum _{k=0}^{n-1}U^kf,\quad n=1,2,\ldots\ .
$$

From the Birkhoff Ergodic Theorem for Markov operators (see
\cite{Fog} for the details) and from the Proposition above we get
the following theorem.

\begin{theorem}\label{Birkhoff}
Suppose $S:(X,\mathcal{B},\mu)\rightarrow (X,\mathcal{B},\mu)$ is
a measure preserving $m$-transformation and $f\in \mathcal{L}^1$.
Then there exists a function $f^*\in \mathcal{L}^1$ such that
$$
A_n(f)\rightarrow f^*, \mu - a.e.
$$
Furthermore, $Uf^*=f^*$ $\mu$-a.e. and $\int _X f^*\ {\rm d}\mu
=\int _X f\ {\rm d}\mu$.
\end{theorem}

\begin{corollary}
Let $1\leq p<\infty$ and let $S$ be a measure preserving
$m$-transformation on $(X,\mathcal{B},\mu)$. If $f\in
\mathcal{L}^p$, then there exists $f^*\in \mathcal{L}^p$ such that
$Uf^*=f^*$ $\mu$-a.e. and $\| f^*-A_n(f)\|_p\rightarrow 0$ as
$n\rightarrow \infty$.
\end{corollary}

$\blacktriangleright $\ Let us fix $1\leq p\leq \infty$ and $f\in
\mathcal{L}^p$. Since $\| A_n(f)\|_p\leq \| f\|_p$, we have by
Fatou's lemma,
$$
\int \limits_X |f^*|^p\ {\rm d}\mu \leq \liminf_{n\rightarrow
+\infty}\int \limits_X |A_n(f)|^p\ {\rm d}\mu \leq \int \limits_X
|f|^p\ {\rm d}\mu \ .
$$
Hence, the operator $L:\mathcal{L}^p\rightarrow \mathcal{L}^p$
defined by $L(f)=f^*$ is a contraction on $\mathcal{L}^p$. By
Theorem \ref{Birkhoff} $\| f^*-A_n(f)\|_p\rightarrow 0$ as
$n\rightarrow \infty$ for any bounded function $f\in
\mathcal{L}^p$. Let $f\in \mathcal{L}^p$ be a function, not
necessarily bounded. For any $\varepsilon >0$ we can find a
bounded function $f_B\in \mathcal{L}^p$ such that $\|
f-f_B\|_p<\varepsilon$. Then, since $L$ is a contraction on
$\mathcal{L}^p$, we have
$$
\| f^*-A_n(f)\|_p\leq \| f_B^*-A_n(f_B)\|_p+\| A_n(f-f_B)\|_p+\|
(f-f_B)^*\|_p\ ,
$$
which can be made arbitrarily small.\ $\blacktriangleleft$

\section{Ergodicity}

Assume $Uf=f$ for some measurable function $f$. It is very
important to know condition on $S$ under that $f$ is constant.

\begin{definition}
We call a nonsingular $m$-transformation $S$ {\bf ergodic} if for
any $B\in \mathcal{B}$, such that $B\backslash
S^{-1}(B)=B^c\backslash S^{-1}(B^c)=\emptyset$, $\mu (B)=0$ or
$\mu (B^c)=0$.
\end{definition}

It is obvious that if $S$ is the union of $\mu$-preserving
$m$-transformations (see Proposition \ref{union}) one of which is
not ergodic, then $S$ is not ergodic.

\begin{theorem}\label{equivalent}
The following three statements are equivalent for any nonsingular
$m$-transformation $S:X\rightarrow X$.
\begin{enumerate}
\item $S$ is ergodic
\item for any $B\in \mathcal{B}$, such that $\mu (B\backslash S^{-1}(B))=\mu
(B^c\backslash S^{-1}(B^c))=0$, $\mu (B)=0$ or $\mu (B^c)=0$.
\item for any disjoint sets $B_1, B_2\in \mathcal{B}$, such that $\mu (B_1\backslash S^{-1}(B_1))=\mu
(B_2\backslash S^{-1}(B_2))=0$, $\mu (B_1)=0$ or $\mu (B_2)=0$.
\end{enumerate}
\end{theorem}

$\blacktriangleright $\ We see at once that $(3)\!\!\Rightarrow \!\! (1)$.

$(1)\!\!\Rightarrow \!\! (2)$ Suppose $S$ is ergodic and $B\in
\mathcal{B}$, such that $\mu (B\backslash S^{-1}(B))=\mu
(B^c\backslash S^{-1}(B^c))=0$. Let $A_1=(B\cap S^{-1}(B))\cup
(B^c\backslash S^{-1}(B^c)),\ A_i=A_{i-1}\cap S^{-1}(A_{i-1})$ for
$i\geq 2$, and $A=\cap _{i=1}^\infty A_i$. We have $A_1\supset
A_2\supset \ldots$ and
$$
A_{i-1}\backslash A_i\subset S^{-1}(A_{i-2}\backslash
A_{i-1})\subset \ldots \subset S^{-i+2}(A_{1}\backslash
A_{2})\subset S^{-i+1}(B\backslash S^{-1}(B))\ .
$$
Therefore, $\mu (A\triangle B)=0$. Let $x\in A$, then there is at
least one point in $S(x)$ that belongs to infinite many of $A_i$.
This gives $A\subset S^{-1}(A)$.

Let $C_1=A^c,\ C_i=C_{i-1}\cap S^{-1}(C_{i-1})$ for $i\geq 2$, and
$C=\cap _{i=1}^\infty C_i$. We have $C_1\supset C_2\supset \ldots$
and
$$
C_{i-1}\backslash C_i\subset \ldots \subset
S^{-i+2}(C_{1}\backslash C_{2})\subset S^{-i+1}(B^c\backslash
S^{-1}(B^c))\cup S^{-i+2}(B\backslash A)\ .
$$
Therefore, $\mu (C\triangle B^c)=0$. Let $x\in C$, then there is
at least one point in $S(x)$ that belongs to infinite many of
$C_i$. This gives $C\subset S^{-1}(C)$. Moreover,
$$
C^c=A\cup C_1\backslash C\subset S^{-1}(A)\cup
S^{-1}(C_1\backslash C)\cup S^{-1}(A)=S^{-1}(C^c)\ .
$$
We conclude from the ergodicity of $S$ that $\mu (B^c)=\mu (C)=0$
or $\mu (B)=\mu (C^c)=0$.

$(2)\!\!\Rightarrow \!\! (3)$ Suppose $(2)$ holds true and let
$B_1, B_2\in \mathcal{B}$ be the disjoint sets, such that $\mu
(B_1\backslash S^{-1}(B_1))=\mu (B_2\backslash S^{-1}(B_2))=0$.
Let $C_1=B_1^c,\ C_i=C_{i-1}\cap S^{-1}(C_{i-1})$ for $i\geq 2$,
and $C=\cap _{i=1}^\infty C_i$. We have $C_1\supset C_2\supset
\ldots$ and $\mu (B_2\backslash C_i)=0$. Therefore $\mu (C)\geq
\mu (B_2)$. Let $x\in C$, then there is at least one point in
$S(x)$ that belongs to infinite many of $C_i$. This gives
$C\subset S^{-1}(C)$. Moreover $\mu (C^c\backslash S^{-1}(C^c))=0$
and $\mu (C^c)\geq \mu (B_1)$. By assumption $\mu(C)=0$ or
$\mu(C^c)=0$. This finishes the proof.\ $\blacktriangleleft$

\begin{example}\label{ex5}
We will prove the ergodisity of
$$
S(x)=\left\{
\begin{array}{cc}
  \{ 2x,1-2x \}\ , & x\in [0,\frac{1}{2}]\phantom{1}  \\
   & \\
  \{ 2x-1 \}\ , & x\in (\frac{1}{2},1]\ .
\end{array}
\right.
$$
\end{example}

$\blacktriangleright $\ Let
\begin{equation}\label{ergex}
B\subset S^{-1}(B)\ \mbox{and}\ B^c\subset S^{-1}(B^c)\ .
\end{equation}
Set $A_1=\{ x:\{ x,1-x \}\subset B\}$, $A_2=\{ x:\{ x,1-x
\}\subset B^c\}$ and $A_3=(A_1\cup A_2)^c$.

Let $x\in A_1$. By (\ref{ergex})
$$
\frac{1+x}{2}\in B \ ,\ \frac{2-x}{2}\in B \ ,\ \frac{1-x}{2}\in B
\ ,\ \frac{x}{2}\in B \ .
$$
Therefore $\bar{S}^{-1}(A_1)\subset A_1$, where $\bar{S}$ is  the
well known ergodic single-valued transformation $\bar{S}(x)=2x\
\mbox{(mod 1)}$, $x\in [0,1]$. By ergodicity of $\bar{S}$,
$\lambda (A_1)=0$ or $\lambda (A_1)=1$. Similarly, $\lambda
(A_2)=0$ or $\lambda (A_2)=1$.

Since $\lambda (A_1)=1$ leads to $\lambda (B^c)=0$ and $\lambda
(A_2)=1$ leads to $\lambda (B)=0$, we need only consider
\begin{equation}\label{ergex2}
\lambda (A_3)=1\ .
\end{equation}
Let $x\in B$. By (\ref{ergex}) and (\ref{ergex2})
$$
\frac{1+x}{2}\in B \ ,\ \frac{2-x}{2}\in B^c \ \mbox{a.s.},\
\frac{1-x}{2}\in B^c \ \mbox{a.s.},\ \frac{x}{2}\in B \
\mbox{a.s.}
$$
Therefore $\lambda(\bar{S}^{-1}(B)\backslash B)$. By ergodicity of
$\bar{S}$, $\lambda (B)=0$ or $\lambda (B)=1$.\
$\blacktriangleleft$

\begin{example}
The 2-transformation $S:[0,1]\rightarrow [0,1]$
$$
S(x)=\left\{
\begin{array}{cc}
  \{ \frac{3}{2}x \}\ , & x\in [0,\frac{1}{3})\phantom{1}  \\
   & \\
  \{ \frac{3}{2}x, \frac{3}{2}x-\frac{1}{2} \}\ , & x\in
  [\frac{1}{3},\frac{2}{3}] \\ & \\
  \{ \frac{3}{2}x-\frac{1}{2} \}\ , & x\in (\frac{2}{3},1]\ .
\end{array}
\right.
$$
is not ergodic.
\end{example}

$\blacktriangleright $\ For instance, $ [0,\frac{1}{2})\subset
S^{-1}([0,\frac{1}{2}))\ \mbox{and}\ [\frac{1}{2},1]\subset
S^{-1}([\frac{1}{2},1])\ .\ \blacktriangleleft $

\begin{proposition}\label{Uinvariant}
Let $S$ be ergodic. If $f$ is measurable and $(Uf)(x)=f(x)$ a.e.,
then $f$ is constant a.e.
\end{proposition}

$\blacktriangleright $\ For each $r\in \mathbb{R}$, $E_r=\{ x\in
X: (Uf)(x)=f(x)>r \}$ is measurable. Then $E_r\subset S^{-1}(E_r)$
and $E_r^c\subset S^{-1}(E_r^c)$, hence $E_r$ has measure 0 or 1.
But if $f$ is not constant a.e., there exists an $r\in \mathbb{R}$
such that $0<\mu (E_r)<1$. Therefore $f$ must be constant a.e.\
$\blacktriangleleft $

\begin{corollary}\label{const}
If a measure preserving $m$-transformation $S$ is ergodic and
$f\in \mathcal{L}^1$, then the limit of the averages $f^*=\int_X
f\ {\rm d}\mu$ is constant a.e. Thus, if $\mu (B)>0$, then for
almost all $x\in X$ there is a orbit of $x$ that returns
infinitely often to $B$.
\end{corollary}

$\blacktriangleright $\ We conclude from Theorem \ref{Birkhoff}
and from Proposition \ref{Uinvariant}, that $ f^*=\int_X f\ {\rm
d}\mu $. To prove the second statement we consider $f=\chi_B$ and
apply Corollary \ref{Poincare}. \ $\blacktriangleleft $

\begin{corollary}\label{uncount}
Let measure preserving $m$-transformation $S$ be ergodic and $\mu
(S_{11}^{-1}(X))<1$, i.e., the set $\{x\in X: |S(x)|\geq 2 \}$ has
positive measure. If $\mu (B)>0$, then for almost all $x\in X$
there are uncountable many orbits of $x$ that return infinitely
often to $B$.
\end{corollary}

$\blacktriangleright $\ We just apply the corollary above to the
sets $B$ and $(S_{11}^{-1}(X)^c$.\ $\blacktriangleleft $

\begin{corollary}\label{leq}
Let $S$ be a measure preserving ergodic $m$-transformation and
$f\in \mathcal{L}^1$ such that $f(x)\geq f(y) (f(x)\leq f(y))$,
for any $y\in S(x)$. Then $f$ is constant a.e.
\end{corollary}
$\blacktriangleright $\ We have $Uf\leq f$, hence the limit of
averages $f^*\leq f$. By Corollary \ref{const} $f=f^*$ is constant
a.e. \ $\blacktriangleleft $

\section{The Frobenius-Perron operator}

Assume that a nonsingular $m$-transformation $S:X\rightarrow X$ on
a normalized measure space is given. We define an operator
$P:\mathcal{L}^1\rightarrow \mathcal{L}^1$ in two steps.

1. Let $f\in \mathcal{L}^1$ and $f\geq 0$. Write
$$
\nu (B)=\int \limits_X f(x)K(x,B)\ {\rm d}\mu\ .
$$
Then, by the Radon-Nikodym Theorem, there exists a unique element
in $\mathcal{L}^1$, which we denoted by $Pf$, such that
$$
\nu (B)=\int \limits_B Pf\ {\rm d}\mu\ .
$$

2. Now let $f\in \mathcal{L}^1$ be arbitrary, not necessarily
nonnegative. Write $f=f^+-f^-$ and define $Pf=Pf^+-Pf^-$. From
this definition we have
$$
\int \limits_B Pf\ {\rm d}\mu =\int \limits_X f^+(x)K(x,B)\ {\rm
d}\mu -\int \limits_X f^-(x)K(x,B)\ {\rm d}\mu
$$
or, more completely,
\begin{equation}\label{Frob}
\int \limits_B Pf\ {\rm d}\mu =\int \limits_X f(x)K(x,B)\ {\rm
d}\mu\ .
\end{equation}

\begin{definition}
If $S:X\rightarrow X$ is a nonsingular $m$-transformation the
unique operator $P:\mathcal{L}^1\rightarrow \mathcal{L}^1$ defined
by equation (\ref{Frob}) is called the {\bf Frobenius-Perron
operator} corresponding to $S$.
\end{definition}
It is straightforward to show that $P$ is a positive linear
operator and
$$
\int \limits_X Pf\ {\rm d}\mu =\int \limits_X f\ {\rm d}\mu\ .
$$

\begin{proposition}
If $f\in \mathcal{L}^1$ and $g\in \mathcal{L}^\infty$, then
$\langle Pf,g\rangle =\langle f,Ug\rangle$, i.e.,
\begin{equation}\label{adjoint}
\int \limits_X (Pf)\cdot g\ {\rm d}\mu =\int \limits_X f\cdot
(Ug)\ {\rm d}\mu\ .
\end{equation}
\end{proposition}

$\blacktriangleright $\ Let $B$ be a measurable subset of $X$ and
$g=\chi _B$. Then the left hand side of (\ref{adjoint}) is
$$
\int \limits_B Pf\ {\rm d}\mu =\int \limits_X f(x)K(x,B)\ {\rm
d}\mu
$$
and the right hand side is
$$
\int \limits_X f\cdot (U\chi _B)\ {\rm d}\mu =\int \limits_X
f\cdot\left( \int \limits_X \chi _B\ {\rm d}K(x,\cdot)\ \right)\
{\rm d}\mu\ =\int \limits_X f(x)K(x,B)\ {\rm d}\mu\ .
$$
Hence (\ref{adjoint}) is verified for characteristic functions.
Since the linear combinations of characteristic functions are
dense in $\mathcal{L}^\infty $, (\ref{adjoint}) holds for all
$f\in \mathcal{L}^1$ and $g\in \mathcal{L}^\infty$.\
$\blacktriangleleft $

The following proposition says that a density $f_*$ is a fixed
point of $P$ if and only if  it is a density of a $S$-invariant
measure $\nu$, absolutely continuous with respect to a measure
$\mu$.

\begin{proposition}\label{newmeasure}
Let $S:X\rightarrow X$ be nonsingular and let $f_*\in
\mathcal{L}^1$ be a density function on $(X,\mathcal{B},\mu)$.
Then $Pf_*=f_*$ a.e., if and only if the measure $\nu =f_*\cdot
\mu$, defined by $\nu (B)=\int _B f_*\ {\rm d}\mu $, is
$S$-invariant.
\end{proposition}

$\blacktriangleright $\ Let $B\subset X$ be a measurable. Then
$$
S\nu (B)=\int \limits_X K(x,B)\ {\rm d}\nu =\int \limits_X
f_*(x)K(x,B)\ {\rm d}\mu =\int \limits_B Pf_*\ {\rm d}\mu\ .
$$
On the other hand
$$
\nu (B)=\int \limits_B f_*\ {\rm d}\mu \ .\ \blacktriangleleft
$$

\begin{proposition}\label{supp}
Let $S:X\rightarrow X$ be a nonsingular $m$-transformation and $P$
the associated Frobenius-Perron operator. Assume that an $f\geq
0,\ f\in \mathcal{L}^1$ is given. Then
$$
{\rm supp}\, f\subset S^{-1}({\rm supp}\, Pf)\ \mbox{a.s.}
$$
\end{proposition}

$\blacktriangleright $\ By the definition of the Frobenius-Perron
operator, we have $Pf(x)=0$ a.e. on $B$ implies that $f(x)=0$ for
a.a. $x\in S^{-1}(B)$. Now setting $B=({\rm supp}\, f)^c$, we have
$Pf(x)=0$ for a.a. $x\in B$ and, consequently, $f(x)=0$ for a.a.
$x\in S^{-1}(B)$, which means that ${\rm supp}\, f\subset
(S^{-1}(B))^c$. Since $(S^{-1}(B))^c\subset S^{-1}(B^c)$ a.s.,
this completes the proof.\ $\blacktriangleleft $

\begin{proposition}
Let $S:X\rightarrow X$ be a nonsingular $m$-transformation and $P$
the associated Frobenius-Perron operator. If $S$ is ergodic, then
there is at most one stationary density $f_*$ of $P$.
\end{proposition}

$\blacktriangleright $\ Assume that $S$ is ergodic and that $f_1$
and $f_2$ are different stationary densities of $P$. Set
$g=f_1-f_2$, so that $Pg=g$. Since $P$ is a Markov operator, $g^+$
and $g^-$ are both stationary densities of $P$. By assumption,
$f_1$ and $f_2$ are not only different but are also densities we
have $g^+\not\equiv 0$ and $g^-\not\equiv 0$. Set
$$
B_1={\rm supp}\, g^+ \quad \mbox{and}\quad B_2={\rm supp}\, g^-\ .
$$
It is evident that $B_1$ and $B_2$ are disjoint sets and both have
positive measure. By Proposition \ref{supp}, we have
$$
B_1\subset S^{-1}(B_1)\ \mbox{a.s.} \quad \mbox{and}\quad
B_2\subset S^{-1}(B_2)\ \mbox{a.s.}
$$
But, from Theorem \ref{equivalent} it follows that $\mu (B_1)=0$
or $\mu (B_2)=0$.\ $\blacktriangleleft $

\section{Applications and generalization}

We now apply the method of $m$-transformation to the intersection
of two middle-$\beta$ Cantor sets (see \cite{Per} and the
references given there).

Let $\alpha \in [\frac{1}{3}, \frac{2}{3}]$ and $\psi _1(x)=\alpha
x$, $\psi _1(x)=\alpha x+1-\alpha$ be a contracting similarity
maps on $I=[0,1]$ endowed with Lebesgue measure $\lambda$. There
is a unique compact set $C_\alpha\subset I$ which satisfies the
set equation
$$
C_\alpha =\psi_1(C_\alpha)\cup \psi_2(C_\alpha)\ .
$$
It is easily checked that $C_\alpha$ is the middle-$\beta$ Cantor
set for $\beta=1-2\alpha$. Let $x\in I$ and $f(x)$ denotes  the
Hausdorff dimension of the set $C_\alpha \cap (C_\alpha +x)$. Let
$B_{ij}=\psi_i(C_\alpha)\cap \psi_j(C_\alpha +x),\ i,j=1,2$. From
the construction of $C_\alpha $ it follows that
$B_{12}=\emptyset$,
$$
{\rm dim_H}B_{11}={\rm dim_H}B_{22}=\left\{
\begin{array}{cc}
  f(\frac{x}{\alpha})\ , & 0\leq x\leq\alpha \\
   & \\
  0\ , & \alpha <x\leq 1
\end{array}
\right.
$$
and
$$
{\rm dim_H}B_{21}=\left\{
\begin{array}{cc}
0\ , & 0\leq x< 1-2\alpha \\ & \\
  f(-\frac{x}{\alpha}+\frac{1}{\alpha}-1)\ , & 1-2\alpha\leq x<1-\alpha \\
   & \\
  f(\frac{x}{\alpha}-\frac{1}{\alpha}+1)\ , & 1-\alpha\leq x\leq 1 \
  .
\end{array}
\right.
$$
Since $C_\alpha \cap (C_\alpha +x)=B_{11}\cup B_{21}\cup B_{22}$,
we have
\begin{equation}\label{max}
f(x)=\max \{ {\rm dim_H}B_{ij}: i,j=1,2\}=\max \{ f(y):y\in
S(x)\}\ ,
\end{equation}
where
$$
S(x)=\left\{
\begin{array}{cc}
  \{ \frac{x}{\alpha} \}\ , & 0\leq x< 1-2\alpha \\ & \\
  \{ \frac{x}{\alpha}, -\frac{x}{\alpha}+\frac{1}{\alpha}-1 \}\ , & 1-2\alpha\leq x\leq\alpha\\ & \\
  \{ -\frac{x}{\alpha}+\frac{1}{\alpha}-1\}\ , & \alpha <x\leq1-\alpha \\ & \\
  \{ \frac{x}{\alpha}-\frac{1}{\alpha}+1\}\ , & 1-\alpha < x\leq
  1
\end{array}
\right.
$$
(compare with Examples \ref{ex2} and \ref{ex5} under $\alpha
=\frac{1}{2}$).

Using Leibniz's rule, we find the Frobenius-Perron operator
corresponding to $S$:
$$
(Pf)(x)=\left\{
\begin{array}{cc}
  \alpha (f(1-\alpha -\alpha x)+f(1-\alpha +\alpha x)+f(\alpha x))\ , & 0\leq x< \frac{1}{\alpha}-2 \\ & \\
  \alpha (f(1-\alpha -\alpha x)+\frac{1}{2}f(1-\alpha +\alpha x)+\frac{1}{2}f(\alpha x))\ , & \frac{1}{\alpha}-2\leq x\leq
  1\ .
 \end{array}
\right.
$$
Assume there exist a stable point $f_*$ of $P$. Then by
Proposition \ref{newmeasure} the measure $\mu =f_*\cdot \lambda$
is $S$-invariant. If in addition $S:(I,\mathcal{B},\mu)\rightarrow
(I,\mathcal{B},\mu)$ is ergodic, then by (\ref{max}) and Corollary
\ref{leq} $f$ is constant $\mu$-a.e. The same method works in case
of the intersection of two arbitrary self-similar sets.

Using $m$-transformations we can develop a new approach to the
self-similar sets with overlaps (see \cite{Bro}, \cite{Nga}). Let
$\psi_1, \ldots ,\psi_m$ be contracting similarity maps on
$\mathbb{R}^n$, and let $X=\cup_{i=1}^m\psi_i(X)$ be an attractor
of the iterated function system. Given normalized measure $\mu$ on
$X$ we consider $m$-transformation of $X$
$$
S(x)=\bigcup_{\{i:x\in \psi_i(X)\}}\psi_i^{-1}(x)\ .
$$
Assume, using the Frobenius-Perron operator corresponding $S$, we
have found $S$-invariant ergodic measure on $X$. This measure
gives us an interesting information about $X$. For instance, if
the conditions of Corollary \ref{uncount} hold true, we see that
a.a. points of $X$ have uncountable many of addresses (see
\cite{Fal} for details).

From these examples we see, that the main problem of the
investigation is to find an $S$-invariant ergodic measure. To
decide this problem we propose a following generalization of an
$m$-transformation.

Given $m$-transformation $S$ on a normalized measure space
$(X,\mathcal{B},\mu)$ we consider a collection of pairs $\{ S_i,
\alpha_i\}_{i=1}^m$, where $S_i:X\rightarrow X$ are the
single-valued measurable transformations such that
$S(x)=\cup_{i=1}^mS_i(x)$ for any $x\in X$, and
$\alpha_i:X\rightarrow [0,1]$ are the measurable functions such
that $\sum_{i=1}^m\alpha_i(x)=1$ for any $x\in X$. Let us consider
the stochastic kernel
$$
K(x,B)=\sum_{i=1}^m\alpha_i(x)\chi_B(S_i(x))
$$
and a new measure on $X$
$$
S\mu (B)\equiv \int \limits_X K(x,B)\ {\rm d}\mu\ .
$$
If we choose $S_i$ and $\alpha_i$ such that $S\mu=\mu$, we can
employ the results of this paper to the measure preserving
transformation $S$.
\vspace{0.2cm}\\
{\bf ACKNOWLEDGMENTS}\\
We are very grateful to Christoph Bandt for suggesting what
measure $S\mu$ is generated by stochastic kernel.

\end{document}